\documentclass[twoside,a4paper,12pt,centertags,reqno]{amsart} 
\usepackage{amsmath,amssymb,verbatim,vmargin}
\usepackage{color}
\usepackage{tikz}
\usepackage{tikz-cd}
\usetikzlibrary{matrix,calc}
\usepackage{hyperref}
\hypersetup{colorlinks,bookmarks=true,linktocpage=true,citecolor=blue,linkcolor=magenta}

\usepackage{eulervm}   

\allowdisplaybreaks 
\usepackage{color}
\usepackage[all]{xy} 

\usepackage{mathtools}
\usepackage{MnSymbol} 

\theoremstyle{plain}
\newtheorem{thm}{Theorem}[section]
\newtheorem{lem}[thm]{Lemma}

\theoremstyle{definition}
\newtheorem{defn}[thm]{Definition}
\newtheorem{rem}[thm]{Remark}

\usepackage{enumerate}
\usepackage{enumitem}

\newcommand{\bRn}{\mathbb{R}^n}

\newcommand{\pd}{\partial}

\newcommand{\fP}{{\mathbf P}}

\def\barint_#1{\mathchoice
            {\mathop{\vrule width 6pt
height 3 pt depth -2.5pt
                    \kern -9.5pt
\intop \kern -4pt}\nolimits_{#1}}%
            {\mathop{\vrule width 5pt height
3 pt depth -2.6pt
                    \kern -6.5pt
\intop \kern -4pt}\nolimits_{#1}}%
            {\mathop{\vrule width 5pt height
3 pt depth -2.6pt
                    \kern -6pt
\intop \kern -4pt}\nolimits_{#1}}%
            {\mathop{\vrule width 5pt height
3 pt depth -2.6pt
          \kern -6pt \intop \kern -4pt}\nolimits_{#1}}}
          
           \def\bariint_#1{\mathchoice
            {\mathop{\vrule width 15pt
height 3 pt depth -2.5pt
                    \kern -15.8pt
\intop \kern -8pt\intop \kern -4pt}\nolimits_{#1}}%
            {\mathop{\vrule width 9pt height
3 pt depth -2.6pt
                    \kern -10.5pt
\intop \kern -8pt\intop \kern -4pt}\nolimits_{#1}}%
            {\mathop{\vrule width 9pt height
3 pt depth -2.6pt
                    \kern -10pt
\intop \kern -8pt\intop \kern -4pt}\nolimits_{#1}}%
            {\mathop{\vrule width 9pt height
3 pt depth -2.6pt
          \kern -8pt \intop \kern -10pt\intop \kern -4pt}
      \nolimits_{  #1}}}

\def\barintlim_#1{\mathchoice
            {\mathop{\vrule width 6pt
height 3 pt depth -2.5pt
                    \kern -8.8pt
\intop \kern -4pt}\limits_{#1}}%
            {\mathop{\vrule width 5pt height
3 pt depth -2.6pt
                    \kern -6.5pt
\intop \kern -4pt}\limits_{#1}}%
            {\mathop{\vrule width 5pt height
3 pt depth -2.6pt
                    \kern -6pt
\intop \kern -4pt}\limits_{#1}}%
            {\mathop{\vrule width 5pt height
3 pt depth -2.6pt
          \kern -6pt \intop \kern -4pt}\limits_{#1}}}
          
           \def\bariintlim_#1{\mathchoice
            {\mathop{\vrule width 15pt
height 3 pt depth -2.5pt
                    \kern -15.8pt
\intop \kern -8pt\intop \kern -4pt}\limits_{#1}}%
            {\mathop{\vrule width 9pt height
3 pt depth -2.6pt
                    \kern -10.5pt
\intop \kern -8pt\intop \kern -4pt}\limits_{#1}}%
            {\mathop{\vrule width 9pt height
3 pt depth -2.6pt
                    \kern -10pt
\intop \kern -8pt\intop \kern -4pt}\limits_{#1}}%
            {\mathop{\vrule width 9pt height
3 pt depth -2.6pt
          \kern -8pt \intop \kern -10pt\intop \kern -4pt}
      \limits_{  #1}}}
          
\renewcommand{\iint}{\int \kern -8pt\int}       

\numberwithin{equation}{section}
\makeatletter
\@namedef{subjclassname@2020}{\textup{2020} Mathematics Subject Classification}
\makeatother

\title{On eigenvalue inequalities of Schmuckenschl\"ager}

\author{Yi C. Huang} 
\address{Universit\'e Sorbonne Paris Nord, Institut Galil\'ee, LAGA, CNRS (UMR 7539), F-93430 Villetaneuse, France} 
\address{School of Mathematical Sciences, Nanjing Normal University, Nanjing 210023, People's Republic of China}
\email{Yi.Huang.Analysis@gmail.com}
\urladdr{https://orcid.org/0000-0002-1297-7674}

\date{\today} 

\thanks{Research of the author is supported by the National NSF grant of China (no. 11801274). 
This paper is completed while the author is on leave, funded by CSC Postdoctoral/Visiting Scholar Program (no. 202006865011), 
at LAGA of Universit\'e Sorbonne Paris Nord.}

\begin{document}

\begin{abstract}
About ten years ago, Schmuckenschl\"ager proved that the lowest eigenvalue of Dirichlet Laplacian for the intersection of two balls 
(i.e., convex, symmetric and compact subsets of $\bRn$ with non-empty interior)
is less than the sum of the lowest eigenvalue for each.
His arguments rely on Kac's formula, the log-concavity of Gaussian measures, the symmetry of balls and Lieb's classical result for the intersection of two domains.
In this note we revisit Schmuckenschl\"ager's proof and propose a direct and elementary Lieb-free approach to these inequalities.
\end{abstract}

\maketitle


\section{Introduction}

Let $\Omega$ be a bounded Lipschitz domain in $\bRn$.
By $\lambda_1(\Omega)$ we denote the first eigenvalue of the following Dirichlet boundary value problem on $\Omega$:
$$-\frac12\Delta u=\lambda_1(\Omega)u\quad\text{and}\quad u|_{\pd\Omega}=0.$$
We are interested in the case when $\Omega$ is the intersection of two balls. 

\begin{defn}
$B\subset\bRn$ is said to be a \textit{ball}, if it is convex, symmetric (more precisely, centrally symmetric about the origin 0), 
compact and has non-empty interior. 
\end{defn}

\begin{thm}[Schmuckenschl\"ager, 2011] \label{thm:SchmEigen}
Let $B_1$ and $B_2$ be two balls in $\bRn$. Then
\begin{equation} \label{eqn:SchmEigen}
\lambda_1(B_1\cap B_2)< \lambda_1(B_1)+\lambda_1(B_2).
\end{equation}
\end{thm}

This subadditivity property can be found in \cite[Proposition 2.5]{Sch11}.
Here we state it as an exact inequality, and this is implicit in Schmuckenschl\"ager's arguments.
Similar results for complex balls under complex interpolation (in the spirit of Cordero-Erausquin \cite{Cor02} and Berndtsson \cite{Ber98}) are also proved there.

Schmuckenschl\"ager's original arguments rely on Kac's formula, the log-concavity of Gaussian measures, the symmetry of balls 
and Lieb's classical result for the intersection of two domains.
In this note we revisit Schmuckenschl\"ager's proof and also propose a direct and elementary Lieb-free approach to these inequalities.

\section{Proof of Theorem \ref{thm:SchmEigen} after Schmuckenschl\"ager}

In formulating the Dirichlet problem, the factor $\frac12$ in front of $\Delta$ is adopted since we shall use a probabilistic approach.
Let $\{X_t, t\geq0\}$ stand for standard Brownian motion in $\bRn$.
For $x\in\bRn$, $\fP^x$ denotes the probability measure such that
$$\fP^x(X_0=x)=1,$$
namely, under $\fP^x$, the process $X_t$ starts at $x$.
Let $T_\Omega$ be the exit time of Brownian motion $X_t$ from $\Omega$.
Kac's formula indicates that for all $x\in\Omega$,
\begin{equation} \label{eqn:Kac}
\lambda_1(\Omega)=-\lim_{t\rightarrow\infty}\frac{\log \fP^x(T_\Omega>t)}{t}.
\end{equation}
One has the following log-concavity result (see also Borell \cite{Bor75}).

\begin{lem}[Schmuckenschl\"ager, 2011] \label{lem:logconc}
Let $A$, $B$ and $C$ be open subsets of $\bRn$, $0\leq\lambda\leq1$ such that
$$(1-\lambda)A+\lambda B\subset C.$$
Then
 \begin{equation} \label{eqn:logconc}
\fP^{(1-\lambda)x+\lambda y}(T_C>t)\geq  \fP^x(T_A>t)^{1-\lambda}\fP^y(T_B>t)^{\lambda}.
\end{equation}
\end{lem}

\begin{rem}
In stating this lemma we also included the endpoint cases $\lambda=0$ and $\lambda=1$. 
They simply follow from the monotonicity of the exit time with respect to the set. 
Indeed, if $x\in\Omega_1\subset\Omega_2$, then for any $t>0$,
$$T_{\Omega_1}\leq T_{\Omega_2} \quad\text{and}\quad\fP^x(T_{\Omega_1}>t)\leq\fP^x(T_{\Omega_2}>t).$$
As a consequence, and using Kac's formula, we have the monotonicity
 \begin{equation} \label{eqn:mono}
\Omega_1\subset\Omega_2\Longrightarrow\lambda_1(\Omega_1)\geq\lambda_1(\Omega_2).
\end{equation}
Moreover, the strict inclusion implies the strict monotonicity.
\end{rem}

\begin{rem}
Complementing \eqref{eqn:SchmEigen}, we have 
$$\lambda_1(B_1\cap B_2)\geq \max\{\lambda_1(B_1),\lambda_1(B_2)\}.$$
\end{rem}

Now we reproduce Schmuckenschl\"ager's proof.
By the symmetry of the balls $B_1$ and $B_2$ (which are also both centred at $0$), 
we have: for any $y\in\bRn$ and any $t>0$,
\begin{equation} \label{eqn:symm}
\fP^0(T_{B_1\cap (B_2+y)}>t)=\fP^0(T_{B_1\cap (B_2-y)}>t),
\end{equation}
hence using Kac's formula, 
\begin{equation} \label{eqn:symm'}
\lambda_1(B_1\cap (B_2+y))=\lambda_1(B_1\cap (B_2-y)).
\end{equation}
Then, by a classical result of Lieb in \cite{Lie83}, 
there exists\footnotemark
\footnotetext{Lieb's result describes a significant intersection between $B_1$ and $B_2+x$.
Therefore, $x$ is small compared to the size of $B_2$.
This observation motivates our next approach.} an $x\in\bRn$ so that 
\begin{equation} \label{eqn:Lieb}
\lambda_1(B_1\cap (B_2+x))< \lambda_1(B_1)+\lambda_1(B_2). 
\end{equation}
Next, by convexity, for any $x\in\bRn$ we have
$$\frac{B_1\cap (B_2+x)+B_1\cap (B_2-x)}{2}\subset B_1\cap B_2,$$
thus by \eqref{eqn:mono}, \eqref{eqn:logconc} plus Kac's formula\footnotemark
\footnotetext{We use $\fP^0$ since $0\in B_1\cap B_2\pm x$, which is reasonable as explained in the previous footnote.}, \eqref{eqn:symm'} and \eqref{eqn:Lieb},
$$\begin{aligned}
\lambda_1(B_1\cap B_2)&\leq \lambda_1\left(\frac{B_1\cap (B_2+x)+B_1\cap (B_2-x)}{2}\right)\\
&\leq \frac{\lambda_1(B_1\cap (B_2+x))+\lambda_1(B_1\cap (B_2-x))}{2}\\
&=\lambda_1(B_1\cap (B_2+x))<\lambda_1(B_1)+\lambda_1(B_2).
\end{aligned}$$
This proves Theorem \ref{thm:SchmEigen}.

\begin{rem}
For recent interests towards Lieb's eigenvalue estimate \eqref{eqn:Lieb}, 
see Frank and Larson \cite{FraLar21} in connection with Davies' Hardy inequality \cite{Dav84}.
\end{rem}

\section{A Lieb-free Proof of Theorem \ref{thm:SchmEigen}}

\textbf{Observation}. We first give the following heuristics:
Lieb's result works for domain intersection in a general context, and
the resulted translation $x\in\bRn$ somehow carrying the information that $B_1$ and $B_2+x$ have the maximal possible intersection.
For this reason it is natural that for the ball intersection we can simply take $x=0$.
In other words, the use of Lieb's intersection result in last proof could be artificial.

\bigskip

\textit{More precisely, our aim is to bypass the use of \eqref{eqn:Lieb} and prove directly}
\begin{equation} \label{eqn:Lieb'}
\lambda_1(B_1\cap B_2)< \lambda_1(B_1)+\lambda_1(B_2),
\end{equation}
\textit{where $B_1$ and $B_2$ are two balls in $\bRn$.}

\bigskip

First, note that the cases $B_1\subset B_2$ and $B_2\subset B_1$ are trivial for \eqref{eqn:Lieb'} and thus can be excluded.
Without loss of generality, we can assume that the aperture $\alpha$ of the cones for each component of $B_1\backslash B_2$ is no greater than that of $B_2\backslash B_1$.
Hence $\alpha\leq\pi/2$.

By homogeneity, we have
$$\frac12\lambda_1(B_1\cap B_2)=\lambda_1(\sqrt2(B_1\cap B_2)).$$
Thus it suffices to prove the following geometric property
\begin{equation} \label{eqn:Geom}
(B_1\backslash B_2)\subset \sqrt2(B_1\cap B_2),
\end{equation}
which implies by convexity the (strict) inclusion
$$\frac{B_1+B_2}{2}\subset \sqrt2(B_1\cap B_2).$$
Indeed, by \eqref{eqn:Geom} and \eqref{eqn:mono} we have
$$\begin{aligned}
\frac12\lambda_1(B_1\cap B_2)&=\lambda_1(\sqrt2(B_1\cap B_2))\\
&< \lambda_1\left(\frac{B_1+B_2}{2}\right)\leq\frac{\lambda_1(B_1)+\lambda_1(B_2)}{2}.
\end{aligned}$$
In the last step, we use \eqref{eqn:logconc} plus Kac's formula as before.

Let us now justify \eqref{eqn:Geom} and Theorem \ref{thm:SchmEigen} then follows.
By the symmetry of $B_1$ and $B_2$, 
it suffices to consider $x\in B_1\backslash B_2$ being the axis of this component.
Denote $$x'=\{\lambda x: \lambda>0\}\cap \pd B_2.$$
Denote by $\rho_1$ the metric that defines $B_1$.
Using $\alpha\leq\pi/2$, together with the symmetry of $B_1$ and $B_2$ again, 
there are $x_1^\sharp, x_2^\sharp \in \pd B_1\cap \pd B_2$ with $|x_1^\sharp|_{\rho_1}=|x_2^\sharp|_{\rho_1}$ so that 
$$\sqrt 2 |x'|_{\rho_1}\geq\sqrt 2\left|\frac{x_1^\sharp+x_2^\sharp}{2}\right|_{\rho_1}\geq |x_1^\sharp|_{\rho_1}.$$
But $|x|_{\rho_1}<|x_1^\sharp|_{\rho_1}$,
thus $x\in \sqrt 2(B_1\cap B_2)$ and \eqref{eqn:Geom} is proved. 

\bibliographystyle{alpha}
 
\bibliography{Hua-DirEigenBalls}

\begin{thebibliography}{Dav84}

\bibitem[Ber98]{Ber98}
Bo~Berndtsson.
\newblock {P}r\'ekopa's theorem and {K}iselman's minimum principle for
  plurisubharmonic functions.
\newblock {\em Mathematische Annalen}, 312(4):785--792, 1998.

\bibitem[Bor75]{Bor75}
Christer Borell.
\newblock The {B}runn-{M}inkowski inequality in {G}auss space.
\newblock {\em Inventiones Mathematicae}, 30(2):207--216, 1975.

\bibitem[CE02]{Cor02}
Dario Cordero-Erausquin.
\newblock {S}antalo's inequality on {$\mathbb{C}^n$} by complex interpolation.
\newblock {\em Comptes Rendus Mathematique}, 334(9):767--772, 2002.

\bibitem[Dav84]{Dav84}
E.~Brian Davies.
\newblock Some norm bounds and quadratic form inequalities for
  {S}chr{\"o}dinger operators. {II}.
\newblock {\em Journal of Operator Theory}, 12:177--196, 1984.

\bibitem[FL21]{FraLar21}
Rupert~L. Frank and Simon Larson.
\newblock Two consequences of {D}avies' {H}ardy inequality.
\newblock {\em Functional Analysis and Its Applications}, 55(2):174--177, 2021.

\bibitem[Lie83]{Lie83}
Elliott~H. Lieb.
\newblock On the lowest eigenvalue of the {L}aplacian for the intersection of
  two domains.
\newblock {\em Inventiones Mathematicae}, 74(3):441--448, 1983.

\bibitem[Sch11]{Sch11}
Michael Schmuckenschl{\"a}ger.
\newblock Inequalities for exit times and eigenvalues of balls.
\newblock {\em Potential Analysis}, 35(3):287--300, 2011.

\end{thebibliography}

\end{document}